\documentclass[12pt]{article}

\usepackage{amsfonts}

\newtheorem{theorem}{Theorem}[section]
\newtheorem{lemma}{Lemma}[section]
\newtheorem{corollary}{Corollary}[section]

\makeatletter
\@addtoreset{equation}{section}
\makeatother

\begin{document}

 \sloppy

\begin{center}
{\bf \Large\bf Jacobi matrices: continued fractions, approximation, spectrum.}
\\ \bigskip
{\large Ianovich E.A.}
\\ \bigskip
{\it first published on website http://yanovich.spb.ru in December 2015}
\\ \bigskip
   {\it Saint Petersburg, Russia}
\\ \bigskip
   {\it\small E-mail: eduard@yanovich.spb.ru}
\end{center}

\begin{abstract}
In this work the spectral theory of self-adjoint operator $A$ represented by Jacobi matrix is considered. The approach is based on the continued fraction representation of the resolvent matrix element of $A$. Different criteria of absolute continuity of a spectrum are found. For the analysis of the absolutely continuous spectrum it is used an approximation of $A$ by a sequence of operators $A_n$ with absolutely continuous spectrum on a given interval $[a,b\,]$ which converges to $A$ in a strong sense on a dense set. In the case when $[a,b\,]\subseteq\sigma(A)$ it was found the sufficient condition of absolute continuity of the operator $A$ spectrum on $[a,b\,]$. This condition uses the notion of equi-absolute continuity. It is constructed the system of functions converging to the distribution function of the operator. In the case of the absolutely continuous spectrum, the system of continuous functions converging to the spectral weight of the operator on a given interval is also constructed and was analyzed. The conditions when the derivative of the distribution function of $A$ belongs to the class $C[a,b\,]$ are also obtained.
\end{abstract}

\bigskip

\section{Introduction.}

Jacobi matrix is the tridiagonal matrix of the form
\begin{equation}
\label{Jacobi_Matrix}
  \left(\begin{array}{ccccc} a_0&b_0&0&0&\ldots\\
                      b_0&a_1&b_1&0&\ldots\\
                      0&b_1&a_2&b_2&\ldots\\
                      0&0&b_2&a_3&\ldots\\
                      \vdots&\vdots&\vdots&\vdots&\ddots
  \end{array}\right),
\end{equation}
where all the $a_n$ are real and all the $b_n$ positive.

An operator in separable Hilbert space $H$ can be assosiated with this matrix as follows~\cite{1}. Let $\{e_n\}_0^{\infty}$ be an orthonormal basis in $H$. Let us define basically the operator $A$ on the basis vectors $e_n$, according to the matrix representation~(\ref{Jacobi_Matrix}):
\begin{equation}
\label{Operator's_action}
\begin{array}{l}
  Ae_n=b_{n-1}e_{n-1}+a_ne_n+b_ne_{n+1}\,,\quad n=1,2,3,\ldots\\
  Ae_0=a_0e_0+b_0e_1
\end{array}
\end{equation}

Further, the operator $A$ is defined by linearity on all finite vectors, the set of which is dense in $H$.
(Finite vector has a finite number of non-zero components in the basis $\{e_n\}_0^{\infty}$.) Due to symmetricity of Jacobi matrix, for any two finite vectors $f$ and $g$ we have
$$
(Af,g)=(f,Ag)
$$

Therefore, the operator $A$ is symmetric and permits a closure. This minimal closed operator is a subject of the present consideration and it is naturally to save for it the notation $A$.

The operator $A$ can have deficiency indices $(1,1)$ or $(0,0)$. In the second case having the main interest in the present paper, operator $A$ is self-adjoint. A simple sufficient condition for self-adjointness of $A$ is Carleman's condition
\begin{equation}
\label{Charleman's_condtion}
\sum\limits_{n=0}^\infty\frac{1}{b_n}=+\infty
\end{equation}
In the self-adjoint case the spectrum of $A$ is simple and $e_0$ is generating element. The information about the spectrum of a self-adjoint operator $A$ is contaned in function
\begin{equation}
\label{resolvent}
R(\lambda)=((A-\lambda E)^{-1}e_0,e_0)=\int\limits_{-\infty}^{\infty}\frac{d\sigma(t)}{t-\lambda}\,,\quad \sigma(t)=(E_{t}e_0,e_0)\,,
\end{equation}
defined at $\lambda\not\in\sigma(A)$. ($E_{t}$ is the resolution of the identity of the operator $A$). The distribution function $\sigma(t)$ contains complete information about the spectrum. It has the following representation
$$
\sigma(\lambda)=\sigma_d(\lambda)+\sigma_{ac}(\lambda)+\sigma_{sc}(\lambda)\,,
$$
where $\sigma_d(\lambda)$ is a saltus function defining eigenvalues, $\sigma_{ac}(\lambda)$ and $\sigma_{sc}(\lambda)$ are the absolutely continuous and the singularly continuous parts respectively.

For the function $R(\lambda)$ there exists an afficient approximation algorithm based on the continued fractions. Namely, the following representation is valid
\begin{equation}
\label{c.fr.representation}
R(\lambda)=\frac{1}{\displaystyle a_0-\lambda-\frac{(b_0)^2}
{\displaystyle a_1-\lambda-\displaystyle
  \frac{(b_1)^2}{\displaystyle a_2-\lambda-\frac{(b_2)^2}{\displaystyle
  a_3-\lambda-\ldots}}}}
\end{equation}
The possibility of such representation and also the conditions and character of convergence are defined by the following theorem~\cite{1}

\begin{theorem}{(Hellinger)}
\label{Hellinger}
{ If the operator $A$ is self-adjoint then the continued fraction in~(\ref{c.fr.representation}) uniformly converges in any closed bounded domain of $\lambda$ without common points with real axis to the analytic function defined by formula~(\ref{resolvent}).}
\end{theorem}

The continued fraction convergence here means the existence of the finite limit of the $n$th approximant
$$
  \frac{1}{\displaystyle a_0-\lambda-\frac{b_0^2}{\displaystyle a_1-
            \lambda-\displaystyle
  \frac{b_1^2}{\displaystyle a_2-\lambda-\ldots-\frac{b_{n-2}^2}{
  \displaystyle
  a_{n-1}-\lambda}}}}=-\frac{Q_n(\lambda)}{P_n(\lambda)}\,,
$$
where $P_n(\lambda)$ and $Q_n(\lambda)$ are 1st and 2nd type polynomials respectively. These polynomials form a pair of linearly independent solutions of a second order finite difference equation
\begin{equation}
\label{recurrent_relations}
b_{n-1}\,y_{n-1}+a_n\,y_n+b_n\,y_{n+1}=\lambda\,y_n\,,\quad (n=1,2,3,\ldots)\,,
\end{equation}
with initial conditions
\begin{equation}
\label{initial_conditions}
P_0(\lambda)=1\,,\quad P_1(\lambda)=\frac{\lambda-a_0}{b_0}\,,\quad
Q_0(\lambda)=0\,,\quad Q_1(\lambda)=\frac{1}{b_0}
\end{equation}
The following equality is also valid
\begin{equation}
\label{vronsqian}
P_{n-1}(\lambda)\,Q_n(\lambda)-P_n(\lambda)\,Q_{n-1}(\lambda)=\frac{1}{b_{n-1}}\,,\quad (n=1,2,3,\ldots)
\end{equation}
Thus
$$
R(\lambda)=-\lim\limits_{n\to\infty}\frac{Q_n(\lambda)}{P_n(\lambda)}\,,\quad \mbox{Im}\,\lambda\ne0
$$

\section{Discrete spectrum.}

Using continued fraction expantion of $R(\lambda)$ and continued fractions convergence conditions it is possible to establish the sufficient conditions for discreteness of operator $A$ spectrum.

Let us begin with the following simple assertion
\begin{lemma}
\label{lemma1}
Denote by $K_n(\lambda)$ the infinite part ("tail") of continued fraction~(\ref{c.fr.representation}) starting with $n$th element of sequences $a_n$ and $b_n$
\begin{equation}
\label{fracion_for_K_n}
K_n(\lambda)=\frac{1}{\displaystyle a_n-\lambda-\frac{(b_n)^2}
{\displaystyle a_{n+1}-\lambda-\displaystyle
  \frac{(b_{n+1})^2}{\displaystyle a_{n+2}-\lambda-\frac{(b_{n+2})^2}{\displaystyle
  a_{n+3}-\lambda-\ldots}}}}
\end{equation}
Then
$$
R(\lambda)=-\frac{Q_n(\lambda)+Q_{n-1}(\lambda)\,b_{n-1}K_n(\lambda)}{P_n(\lambda)+P_{n-1}(\lambda)\,b_{n-1}K_n(\lambda)}\,,\quad
(n=1,2,3,\ldots)
$$
\end{lemma}

{\bf Proof.}\, The proof is obtained by induction on $n$. The validity of this formula for $n=1$ follows from initial conditions ~(\ref{initial_conditions}) for polynomials  $P_n$ и $Q_n$. Assuming the validity of formula
$$
R(\lambda)=-\frac{Q_{n-1}+Q_{n-2}\,b_{n-2}K_{n-1}}{P_{n-1}+P_{n-2}\,b_{n-2}K_{n-1}}\,,
$$
let us prove the same one for the value of index greater on unit. Expressing $K_{n-1}$ through $K_n$, we obtain
$$
R(\lambda)=-\frac{\displaystyle Q_{n-1}+Q_{n-2}\,b_{n-2}\,\frac{1}{\displaystyle a_{n-1}-\lambda-(b_{n-1})^2\,K_{n}}}
{\displaystyle P_{n-1}+P_{n-2}\,b_{n-2}\,\frac{1}{\displaystyle a_{n-1}-\lambda-(b_{n-1})^2\,K_{n}}}=
$$
$$
=-\frac{(a_{n-1}-\lambda)Q_{n-1}+Q_{n-2}\,b_{n-2}-Q_{n-1}\,(b_{n-1})^2\,K_{n-1}}
{(a_{n-1}-\lambda)P_{n-1}+P_{n-2}\,b_{n-2}-P_{n-1}\,(b_{n-1})^2\,K_{n-1}}
$$
Using~(\ref{recurrent_relations}), we find
$$
R(\lambda)=-\frac{-Q_n\,b_{n-1}-Q_{n-1}\,(b_{n-1})^2\,K_{n-1}}
{-P_n\,b_{n-1}-P_{n-1}\,(b_{n-1})^2\,K_{n-1}}=
-\frac{Q_n+Q_{n-1}\,b_{n-1}K_n}{P_n+P_{n-1}\,b_{n-1}K_n}
$$
Thus the formula is valid for all natural $n$.\\
-------------------------------------------------------------------------------

Exactly infinite part of the continued fraction (the behavior $a_n$ and $b_n$ at infinity) defines essentialy  the spectrum of $A$. The discrete component only can qualitatively depends on finite number of first terms of sequences $\{a_n\}$ and $\{b_n\}$. But this possible changes can consist in appearance (or disappearance) of a finite number of eigenvalues.

To state that the spectrum of $A$ is discrete we shall find such conditions that for any domain $D$, no matter how large it is, there exists a number $N_D$ such that the function $K_n(\lambda)$ is analytic in $D$ for $n>N_D$ ($N_D$ depends on $D$). We use for that the known convergence theorems of continued fractions, like Pringsheim's theorem.

First result of such type was established probably by Chihara~\cite{16}.
Note that in the work~\cite{5} another close criterion was obtained by different method.
\begin{theorem}
\label{criteri_discretnosti}
Assume that
$$
\lim\limits_{n\to\infty}|a_n|=\infty\,,\quad
\overline{\lim}\,\frac{\textstyle b_{n-1}^2}{|a_{n-1}a_n|}=c<1/4
$$
Then the operator $A$ (or its self-adjoint extension) has a discrete spectrum.
\end{theorem}

{\bf Proof.}\, Transform the continued fraction~(\ref{fracion_for_K_n}) in the expression for $K_n$ into equivalent one
\begin{equation}
\label{transformed_fraction_for_K_n}
K_n(\lambda)=\frac{C_0(\lambda)}{\displaystyle 2+\frac{C_1(\lambda)}
{2+\displaystyle\frac{C_2(\lambda)}{\displaystyle 2+\ldots}}}\,,
\end{equation}
where
\begin{equation}
\label{C_k}
C_0(\lambda)=\frac{2}{a_n-\lambda}\,,\quad
C_k(\lambda)=\frac{-4\,b_{n+k-1}^2}{(a_{n+k-1}-
\lambda)(a_{n+k}-\lambda)}\,,\qquad (k=1,2,\ldots)
\end{equation}
To avoid a bulkiness we will omit an index $n$ in the writing of $C_k(\lambda)$ although they depend on $n$.

A continued fraction in the formula~(\ref{transformed_fraction_for_K_n}) is defined as usualy as a limit of $n$th approximant~\cite{2,3}
$$
K_n(\lambda)=\lim\limits_{k\to\infty}\frac{C_0(\lambda)}{\displaystyle 2+\frac{C_1(\lambda)}
{2\displaystyle+\ldots+\frac{C_k(\lambda)}{\displaystyle 2}}}=\lim\limits_{k\to\infty}\frac{M_k}{N_k}\,,
$$
where $M_k$ and $N_k$ are $n$th numerator and denominator satisfying the following three-term recurrence relations
\begin{equation}
\label{N_k}
  \begin{array}{cc}
  M_k=2\,M_{k-1}+C_k(\lambda)\,M_{k-2}\: ; & M_{-1}\equiv 0\:,\:M_{-2}\equiv 1 \\
  N_k=2\,N_{k-1}+C_k(\lambda)\,N_{k-2}\: ; & N_{-1}\equiv 1\:,\:N_{-2}\equiv 0
  \end{array}
\end{equation}
The $n$th approximant of the continued fraction~(\ref{transformed_fraction_for_K_n}) coincides with $n$th approximant of the continued fraction~(\ref{fracion_for_K_n}). Therefore these fractions converge or diverge simultaneously and have common limit in the case of convergence.

Using the formula~\cite{2,3}
\begin{equation}
\label{formula_for_N/M}
  \frac{M_k}{N_k}=\sum_{i=0}^k(-1)^i\:\frac{C_0\cdot\ldots
  \cdot C_i}{{N_i}N_{i-1}}\,,
\end{equation}
it is possible to present the function $K_n(\lambda)$ as a series
\begin{equation}
\label{series_for_K_n}
K_n(\lambda)=\sum_{i=0}^\infty(-1)^i\:\frac{C_0\cdot\ldots
  \cdot C_i}{{N_i}N_{i-1}}\,,
\end{equation}
Convergence of the continued fraction~(\ref{transformed_fraction_for_K_n}) is equivalent to the convergence of this series.

In the proof of the convergence below we essentially use the idea of Pringsheim's theorem proof given in~\cite{3}.

It follows from the conditions of a theorem and formulas~(\ref{C_k}) that for any bounded domain $D$ in complex plane there exists a number $N_D$ such that all functions $C_k(\lambda)$ are analytic in $D$, and for any $\lambda\in D$ and $n>N_D$ the following inequalities are fulfilled
$$
|C_k(\lambda)|<1\,,\quad (k=0,1,2,\ldots)
$$
Using this, and formulas~(\ref{N_k}), we obtain for $\lambda\in D$ and $n>N_D$
$$
|N_k|\ge2\,|N_{k-1}|-|C_k|\,|N_{k-2}|\ge2\,|N_{k-1}|-|N_{k-2}|
$$
or
$$
|N_k|-|N_{k-1}|\ge |N_{k-1}|-|N_{k-2}|
$$
Hence it follows that
\begin{equation}
\label{ner-vo_for_N_k}
|N_k|-|N_{k-1}|\ge|N_{-1}|-|N_{-2}|=1,\quad (k=1,2,\ldots)
\end{equation}
Therefore for the partial sums of a series~(\ref{series_for_K_n}) we have the estimate
\begin{equation}
\label{estimation_for_N/M}
\left|\sum_{i=0}^{k}\frac{\displaystyle\prod_{j=0}^iC_i}{N_iN_{i-1}}\right|
\le\sum_{i=0}^{k}\frac{1}{|N_iN_{i-1}|}\le
\sum_{i=0}^{k}\frac{|N_i|-|N_{i-1}|}{|N_iN_{i-1}|}=
\sum_{i=0}^{k}\left(\frac{1}{|N_{i-1}|}-\frac{1}{|N_i|}\right)=1-\frac{1}{|N_k|}\le 1
\end{equation}
It follows that the series~(\ref{series_for_K_n}) converges absolutely, and hence, the continued fraction~(\ref{transformed_fraction_for_K_n}) converges for all $\lambda\in D$. Let us prove that the limit function $K_n(\lambda)$ is analytic in $D$. The functions $C_k(\lambda)$ are analytic in $D$. Therefore, from the recurrence relations~(\ref{N_k}) it follows that $M_k$ and $N_k$ are also analytic functions of $\lambda$ in $D$. Besides, as it follows from the inequality~(\ref{ner-vo_for_N_k}), the functions $N_k(\lambda)$ have no zeros in $D$. Hence the function $K_n(\lambda)$ is a limit at any point $\lambda\in D$ of sequence of analytic in $D$ functions
$$
\frac{M_k(\lambda)}{N_k(\lambda)}
$$
It follows from the formula~(\ref{formula_for_N/M}) and estimate~(\ref{estimation_for_N/M}) that this sequence is uniformly bounded in $D$ and consequently it is compact in this domain. Therefore by Vitali's theorem~\cite{4} this sequence converges uniformly in $D$ and by Weierstrass's theorem the limit function $K_n(\lambda)$ is analytic in this domain. From~(\ref{formula_for_N/M}) and~(\ref{estimation_for_N/M}) it follows also that
$$
|K_n(\lambda)|\le 1\,,\quad \lambda\in D
$$
Thus it is proved that for any bounded domain $D$ of a complex plane there exists a number $N_D$ such that for $n>N_D$ the function $K_n(\lambda)$ is analytic in $D$ and is a limit of the continued fraction~(\ref{fracion_for_K_n}) converging in $D$ .
From Lemma~\ref{lemma1} it follows that the function $R(\lambda)$ is represented as a quotient of the two analytic in $D$ functions (polynomials $P_n(\lambda)$ and $Q_n(\lambda)$ are obviously analytic in any bounded part of a complex plane). Therefore the function $R(\lambda)$ is analytic in $D$, except may be a finite number of points at which the denominator is zero. With the exception of this finite number of points (they can absent at all) the continued fraction~(\ref{c.fr.representation}) converges in $D$. If the operator $A$ is self-adjoint, the discreteness of its spectrum follows from the fact that only a finite number of singular points of the function $R(\lambda)$, which are its poles, can are contained in any bounded domain. By virtue of sef-adjointness of the operator $A$ these poles will exist and will be on a real axis being the eigenvalues of $A$. A unique condensation point of the eigenvalues will be at infinity. If the operator $A$ is not self-adjoint, any its self-adjoint extension has a discrete spectrum since the resolvent is compact~\cite{13}. The theorem is proved.\\
-------------------------------------------------------------------------------

\section{Approximation and absolutely continuous spectrum.}

Let us consider the following problem. Let $\{A_n\}_1^\infty$ be a sequence of linear self-adjoint operators in a separable Hilbert space $H$ with absolutely continuous spectrum on some interval $[a,b\,]$ and $\Phi$ be some dense set in $H$. Suppose that on $\Phi$ there exist in a strong sense the limit $\lim\limits_{n\to\infty}A_n$ and the closure of this limit operator is self-adjoint operator A. That is, for any $e\in\Phi$
$$
\lim\limits_{n\to\infty}\|(A_n-A)e\|=0
$$

Suppose also that $[a,b\,]\subseteq\sigma(A)$. There arises the question: what is the condition which guarantees that the spectrum of a limit operator $A$ is also absolutely continuous on $[a,b\,]$ ?

Let $E_\lambda^{(n)}$ and $E_\lambda$ be the resolutions of the identity of the operators $A_n$ and $A$ respectively and let
$$
\sigma_n(\lambda;\,e)=(E_\lambda^{(n)} e, e),\quad e\in H
$$
$$\sigma(\lambda;\,e)=(E_\lambda e, e),\quad e\in H
$$
By assumption all functions $\sigma_n(\lambda;\,e)$ are absolutely continuous on $[a,b\,]$ so that
$$
\sigma_n(\lambda;\, e)=\sigma_n(a;\, e)+\int\limits_{a}^\lambda f_n(t; e)\,dt\,,\quad e\in H\,,\quad \lambda\in[a,b\,],
$$
where $f_n(t;\,e)$ is a sequence of non-negative summable on $[a,b\,]$ functions.

It may to suppose that the convergence of the sequence $A_n$ to the operator $A$ in a strong sense on the set $\Phi$ should leads to the convergence of the sequence $\sigma_n(\lambda; e)$ to the function $\sigma(\lambda; e)$. It turns out it is true.
\begin{lemma}
\label{lemma2}
If $\lambda$ is not an eigenvalue of the operator $A$, then for all $e\in H$
\begin{equation}
\label{sigma_n_to_sigma}
\lim\limits_{n\to\infty}\sigma_n(\lambda; e)=\sigma(\lambda; e)
\end{equation}
\end{lemma}

{\bf Proof.}\, It is known~\cite{6} that the strong convergence of the sequence $A_n$ to $A$ on a dense subset in $H$ leads to the strong convergence of the sequence $E_\lambda^{(n)}$ to $E_\lambda$ for those $\lambda$ which are not the eigenvalues of $A$. That is, if $\lambda$ is not an eigenvalue of the operator $A$, then for all $e\in H$
$$
\lim\limits_{n\to\infty}\|(E_\lambda^{(n)}-E_\lambda)e\|=0
$$
Further, we have
$$
|\sigma_n(\lambda; e)-\sigma(\lambda; e)|=|((E_\lambda^{(n)}-E_\lambda)e, e)|\le\|(E_\lambda^{(n)}-E_\lambda)e\|\,,
$$
whence the statement of the Lemma readily follows.\\
-------------------------------------------------------------------------------

Note that as we suppose that the space $H$ is separable, the set of eigenvalues of the operator $A$ is countable and hence the equality~(\ref{sigma_n_to_sigma}) is fulfilled for almost all $\lambda$.

Let us give now the answer on the question stated above.
\begin{theorem}
\label{criterion_of_absolutely_continuity}
In order that the spectrum of the operator $A$ on the interval $[a,b\,]$ $(\,[a,b\,]\subseteq\sigma(A)\,)$ be absolutely continuous it is sufficient that for any $e\in H$ the sequence of functions $\{\sigma_n(\lambda;\, e)\}$ be equi-absolutely continuous on $[a,b\,]$, or, that is the same, the sequence of functions $\{f_n(t; e)\}$ has equi-absolutely continuous integrals on $[a,b\,]$. Under these conditions the sequence $\{\sigma_n(\lambda;\, e)\}$ uniformly converges to $\sigma(\lambda;\, e)$ on $[a,b\,]$.
\end{theorem}

{\bf Proof.}\, Fix some vector $e\in H$. For the brevity we will omit this vector in the arguments of all functions.

Consider the sequence of functions $\sigma_n(\lambda)$ on the interval $[a,b\,]$.
$$
\sigma_n(\lambda)=\sigma_n(a)+\int\limits_{a}^\lambda f_n(t)\,dt\,,\quad \lambda\in[a,b\,]
$$

Here $f_n(t)\ge0$ on $[a,b\,]$ but the case $f_n(t)\equiv0$ is also possible if for example the corresponding vector $e$ belongs to the discrete subspace of the operator $A_n$. In this case the function $\sigma_n(\lambda)$ is constant on $[a,b\,]$.

By assumption the sequence $\sigma_n(\lambda)$ is equi-absolutely continuous on $[a,b\,]$. Besides, by virtue of its definition, it is uniformly bounded ($0\le\sigma_n(\lambda)\le\|e\|^2\,,\, \lambda\in R$). Hence~\cite{7}, from the sequence $\sigma_n(\lambda)$ one can extract a partial sequence $\sigma_{n_k}(\lambda)$ which uniformly converges on $[a,b\,]$, moreover the limit function $\sigma_1(\lambda)$ is absolutely continues on $[a,b\,]$. So, we obtain
$$
(*)\quad\lim\limits_{n\to\infty}\sigma_n(\lambda)=\sigma(\lambda)\quad\mbox{for almost all $\lambda$ from $[a,b\,]$}
$$
$$
\lim\limits_{k\to\infty}\sigma_{n_k}(\lambda)=\sigma_1(\lambda)\quad\mbox{for all $\lambda$ from $[a,b\,]$}
$$
Comparing this two limit equalities we can note that $\sigma_1(\lambda)=\sigma(\lambda)$ for almost all $\lambda\in[a,b\,]$. Let us show that the first equality (*) is fulfilled everywhere in $[a,b\,]$. Let $\lambda_0$ be a point in which the equality (*) perhaps is not valid. Consider the sequence $\sigma_n(\lambda_0)$. We have
$$
|\sigma_n(\lambda_0)-\sigma_m(\lambda_0)|\le|\sigma_n(\lambda_0)-\sigma_n(\lambda')|+|\sigma_n(\lambda')-\sigma_m(\lambda')|+
|\sigma_m(\lambda')-\sigma_m(\lambda_0)|
$$
Here $\lambda'$ is an arbitrary point of the interval $[a,b\,]$. Since the set of the possible points $\lambda_0$ has measure zero, in any arbitrarily small neighbourhood of the point $\lambda_0$ there exists an infinite number of points in which the sequence $\sigma_n(\lambda)$ converges. Hence one can choose the point $\lambda'$ such that on the one hand the sequence $\sigma_n(\lambda')$ converges and on the other hand the inequality
$$
|\sigma_n(\lambda_0)-\sigma_n(\lambda')|<\epsilon/3
$$
holds for any $\epsilon>0$ and any $n$.

The validity of this inequality for any $n$ is possible by virtue of equi-absolute continuity of the sequence $\sigma_n(\lambda)$.

Fixing such $\lambda'$ and $\epsilon$, one can choose a number $N$ such that the inequality
$$
|\sigma_n(\lambda')-\sigma_m(\lambda')|<\epsilon/3
$$
is valid, provided $n,m>N$. Then for the same $n$ and $m$
$$
|\sigma_n(\lambda_0)-\sigma_m(\lambda_0)|<\epsilon
$$
Since $\epsilon$ is arbitrarily small the sequence $\sigma_n(\lambda_0)$ converges by Cauchy's criterion.

Thus the sequence $\sigma_n(\lambda)$ converges everywhere in $[a,b\,]$ and therefore $\sigma(\lambda)=\sigma_1(\lambda)$ also everywhere in $[a,b\,]$. Hence the function $\sigma(\lambda)$ is absolutely continuous on this interval. Note that the sequence $\sigma_n(\lambda)$ converges to $\sigma(\lambda)$ uniformly on $[a,b\,]$ (it can be proved by contradiction).

Thus we proved that for any $e\in H$ the functions $\sigma(\lambda)\equiv\sigma(\lambda; e)$ are absolutely continuous on $[a,b\,]$. Since $[a,b\,]\subseteq\sigma(A)$ it follows already that the spectrum of the operator $A$ is absolutely continuous on this interval. The theorem is proved.\\
-------------------------------------------------------------------------------

\begin{corollary}
\label{criterion_of_absolutely_continuity_corollary}
If for any $e\in H$ the sequence of functions $\{\sigma_n(\lambda;\, e)\}$ is equi-absolutely continuous on any finite interval of a real axis then the spectrum of the operator $A$ is purely absolutely continuous.
\end{corollary}

{\bf Remark.} In the conditions of the theorem one can replace the requirement $"$for any  $e\in H$$"$ by the requirement $"$for any $e$ from generation subspace of the operator $A$$"$. In particular, if the spectrum of the operator $A$ is simple, as in the case of Jacobi matrix, then it is sufficient to consider the functions $\{\sigma_n(\lambda;\, e_0)\}$ only, where $e_0$ is a generating vector.\\
------------------------------------------------------------------------------

Turn back to the case when the operator $A$ is represented by Jacobi matrix~(\ref{Jacobi_Matrix}) and choose the approximating sequence $A_n$ as follows
\begin{equation}
\label{A_n}
A_n=\left(\begin{array}{cccccccc} a_0&b_0&\ldots&0&0&0&0&\ldots\\
                      b_0&a_1&\ldots&0&0&0&0&\ldots\\
                      \vdots&\vdots&\ddots&\vdots&\vdots&\vdots&\vdots&\vdots\\
                      0&0&\ldots&a_{n-1}&b_{n-1}&0&0&\ldots\\
                      0&0&\ldots&b_{n-1}&a_n&b_n&0&\ldots\\
                      0&0&\ldots&0&b_n&a_n&b_n&\ldots\\
                      0&0&\ldots&0&0&b_n&a_n&\ldots\\
                      \vdots&\vdots&\vdots&\vdots&\vdots&\vdots&\vdots&\ddots
  \end{array}\right)\,,\quad (n=0,1,\ldots)
\end{equation}
The matrix $A_n$ differs from $A$ by that the elements of the diagonals are not changed starting with $n$rh number. If we denote the elements of the diagonal sequences of $A_n$ by $a^{(n)}_k$ and $b^{(n)}_k$, then
$$
a_k^{(n)}=\left[\begin{array}{c} a_k\,,\quad k<n\\
                                 a_n\,,\quad k\ge n
                \end{array}\right.\,,\quad\quad\quad
                b_k^{(n)}=\left[\begin{array}{c} b_k\,,\quad k<n\\
                                 b_n\,,\quad k\ge n
                \end{array}\right.
$$

Such matrices $A_n$ were used in works~\cite{17,18,19}. Note that they were used also in the work~\cite{11}, where the absolutely continuous spectrum of Jacobi matrices also was considered. However this matrices were used there in connection with commutation relations rather than with approximation. 

It is clear that the matrices $A_n$ describe a sequence of bounded ( $D(A_n)=H$ ) self-adjoint operators (Carleman's condition~(\ref{Charleman's_condtion}) is obviously fulfilled) and it is natural to save the same notation for it.

The sequence $A_n$ strongly converges to the operator $A$ on the set of finite vectors since for any finite vector the difference $A-A_n$ vanishes, provided $n$ is sufficiently large. To use the Theorem~(\ref{criterion_of_absolutely_continuity}) we should state the conditions under which the spectrum of $A_n$ is absolutely continuous. Note that the spectrum of the operators $A_n$ and $A$ is simple and all information about the spectrum contains in the functions $\sigma_n(\lambda)=(E_\lambda^{(n)} e_0, e_0)$ and $\sigma(\lambda)=(E_\lambda e_0, e_0)$ (generating vector $e_0$ is the same for $A$ and $A_n$). Below we will state the absolute continuity conditions of the spectrum of $A_n$. We will give also an explicit expression for the functions $f_n(t)=\sigma'_n(t)$. It is possible to do with the help of continued fraction approach.

Denote by $R_n(\lambda)$ the function $R(\lambda)$ for the operator $A_n$ defined by~(\ref{resolvent}).
Applying the Lemma~(\ref{lemma1}) and taking into account the definition of the operators $A_n$, we can present $R_n(\lambda)$ in the form
\begin{equation}
\label{R_n}
R_n(\lambda)=-\frac{Q_n(\lambda)+Q_{n-1}(\lambda)\,b_{n-1}K_n(\lambda)}{P_n(\lambda)+P_{n-1}(\lambda)\,b_{n-1}K_n(\lambda)}\,,
\end{equation}
where
\begin{equation}
\label{hvost}
K_n(\lambda)=\frac{1}{\displaystyle a_n-\lambda-\frac{(b_n)^2}
{\displaystyle a_{n}-\lambda-\displaystyle
  \frac{(b_{n})^2}{\displaystyle a_{n}-\lambda-\frac{(b_{n})^2}{\displaystyle
  a_{n}-\lambda-\ldots}}}}
\end{equation}

Continued fraction for $K_n(\lambda)$ describes an infinite part of the matrix of the operators $A_n$ in which the elements of the sequences are constant and equal $a_n$ and $b_n$. The function $K_n(\lambda)$ can be regarded as a matrix element of the resolvent of the operator represented by Jacobi matrix $J_n$ with constant elements along the diagonals
$$
J_n=\left(\begin{array}{ccccc} a_n&b_n&0&0&\ldots\\
                      b_n&a_n&b_n&0&\ldots\\
                      0&b_n&a_n&b_n&\ldots\\
                      0&0&b_n&a_n&\ldots\\
                      \vdots&\vdots&\vdots&\vdots&\ddots
  \end{array}\right)
$$
This Jacobi matrix is naturally connected with Tchebychev's polynomials, and the function $K_n(\lambda)$ can be found explicitly.
Note first that the operator $J_n$ is self-adjoint and by the theorem~(\ref{Hellinger}) the continued fraction~(\ref{hvost}) converges  to the function $K_n(\lambda)$ for any non-real $\lambda$. Write next the identity for the $k$th and $(k-1)$th approximants which follows readily from the structure of the continued fraction~(\ref{hvost})
$$
K_n^{(k)}(\lambda)=\frac{1}{\displaystyle a_n-\lambda-(b_n)^2\,K_n^{(k-1)}(\lambda)}
$$
Letting here $\mbox{Im}\,\lambda\ne0$ and passing to the limit, $k\to\infty$, one obtains the simple equation for $K_n(\lambda)$
$$
K_n(\lambda)=\frac{1}{\displaystyle a_n-\lambda-(b_n)^2\,K_n(\lambda)}\,,
$$
which has a solution
$$
K_n(\lambda)=\frac{a_n-\lambda+\sqrt{(a_n-\lambda)^2-4(b_n)^2}}{2\,(b_n)^2}
$$
One of the two branches of the square root is chosen from additional condition
$$
\mbox{Im}\,K_n(\lambda)>0\,,\quad\mbox{Im}\,\lambda>0\,,
$$
which is valid for a matrix element of the resolvent of any self-adjoint operator~\cite{6}. Using this one obtains for any $x\in\mathbb R$
\begin{equation}
\label{limit_K_n1}
\lim\limits_{\epsilon\to+0}K_n(x+i\epsilon)=D_n(x)+i\,B_n(x)
\end{equation}
\begin{equation}
\label{limit_K_n2}
D_n(x)=\left[\begin{array}{cr}\frac{\displaystyle a_n-x}{\displaystyle 2\,(b_n)^2}\,,& |a_n-x|\le2b_n\\
                             \frac{\displaystyle a_n-x\pm\sqrt{(a_n-x)^2-4(b_n)^2}}{\displaystyle 2\,(b_n)^2}\,,& |a_n-x|\ge2b_n
             \end{array}\right.
\end{equation}
\begin{equation}
\label{limit_K_n3}
B_n(x)=\left[\begin{array}{cr}\frac{\displaystyle\sqrt{4(b_n)^2-(a_n-x)^2}}{\displaystyle 2\,(b_n)^2}\,,& |a_n-x|\le2b_n\\
                             0\,,& |a_n-x|\ge2b_n
             \end{array}\right.
\end{equation}
In the expression for $D_n(x)$ the sign $"$+$"$ before the square root is taken when $x-a_n>2b_n$, the sing $"$$-$$"$ when $a_n-x>2b_n$.

From~(\ref{limit_K_n1}) - (\ref{limit_K_n3}) it follows that the spectrum of $J_n$ is absolutely continuous and concentrated on the interval $\sigma(J_n)=[a_n-2b_n,a_n+2b_n]$. The function $K_n(\lambda)$ is the matrix element of the resolvent of self-adjoint operator $J_n$. Hence for the function $K_n(\lambda)$ one has the integral representation
$$
K_n(\lambda)=\int\limits_{-\infty}^{\infty}\frac{f_{J_n}(x)\,dx}{x-\lambda}=
\frac{1}{2\pi (b_n)^2}\int\limits_{a_n-2b_n}^{a_n+2b_n}\frac{\sqrt{4(b_n)^2-(a_n-x)^2}\,dx}{x-\lambda}\,,
$$
where
\begin{equation}
\label{f_J_n}
f_{J_n}(x)=\left[\begin{array}{cr}\frac{\displaystyle\sqrt{4(b_n)^2-(a_n-x)^2}}{\displaystyle 2\pi\,(b_n)^2}\,,& |a_n-x|\le2b_n\\
                             0\,,& |a_n-x|>2b_n
             \end{array}\right.
\end{equation}
The function $f_{J_n}(x)$ at $a_n=0$, $b_n=1/2$ is the spectral weight for Tchebychev's polynomials~\cite{8,9}. The element $a_n$ defines the center of the spectrum and $b_n$ defines the width of the spectrum. The center and the width of the spectrum change according to the behavior of $a_n$ and $b_n$.

Prove that under some conditions the spectrum of $A_n$ is the same as the spectrum of $J_n$.

\begin{theorem}
\label{spectr_A_n}
 Inside the interval $I_n=[a_n-2b_n,a_n+2b_n]$ the spectrum of the operator $A_n$ is always absolutely continuous. The spectral weight $f_n(x)$ with accuracy to the values on a set of measure zero is defined by
\begin{equation}
\label{f_n(x)_3.6}
f_n(x)=\frac{\frac{\displaystyle 1}{\displaystyle \pi}\,B_n(x)}{P^2_n(x)-\frac{\displaystyle b_{n-1}}{\displaystyle b_n}P_{n-1}(x)P_{n+1}(x)}=
\frac{f_{J_n}(x)}{P^2_n(x)-\frac{\displaystyle b_{n-1}}{\displaystyle b_n}P_{n-1}(x)P_{n+1}(x)}\,,
\end{equation}
where $f_{J_n}(x)$ is spectral weight of $J_n$.

If for any $n\in\mathbb N$ the inequalities
$$
b_{n}\ge b_{n-1}\,,\quad |a_n-a_{n-1}|\le2(b_n-b_{n-1})
$$
hold, then the spectrum of $A_n$ is purely absolutely continuous, concentrated on the interval $I_n$ and $I_n\subseteq I_{n+1}$.
\end{theorem}

{\bf Proof.}\, Prove first that the absolutely continuous part of the spectrum of $A_n$ coincides with the spectrum of
$J_n$. For this purpose transform the expression~(\ref{R_n}) for $R_n(\lambda)$
$$
R_n(\lambda)=-\frac{Q_{n-1}(\lambda)}{P_{n-1}(\lambda)}-\frac{Q_{n}(\lambda)\,P_{n-1}(\lambda)-Q_{n-1}(\lambda)\,P_{n}(\lambda)}
{P_{n-1}(\lambda)\,(P_n(\lambda)+P_{n-1}(\lambda)\,b_{n-1}K_n(\lambda))}
$$
or, using~(\ref{vronsqian}),
$$
R_n(\lambda)=-\frac{Q_{n-1}(\lambda)}{P_{n-1}(\lambda)}-\frac{1}{b_{n-1}\,P_{n-1}(\lambda)\,
(P_n(\lambda)+P_{n-1}(\lambda)\,b_{n-1}K_n(\lambda))}
$$
Let $\lambda=x+i\epsilon\,,\:x\in\mathbb R$. Using~(\ref{limit_K_n1}), (\ref{limit_K_n3}) one obtains
\begin{equation}
\label{limit_R_n}
\lim\limits_{\epsilon\to+0}\mbox{Im}\,R_n(x+i\epsilon)=\frac{B_n(x)}
{(P_n(x)+b_{n-1}\,P_{n-1}(x)\,D_n(x))^2+(b_{n-1}\,P_{n-1}(x)\,B_n(x))^2}
\end{equation}
It is easy to note that the denominator of the fraction~(\ref{limit_R_n}) nowhere within the interval of the absolutely continuous spectrum of $J_n$ is zero. Indeed, if $|x-a_n|<2b_n$, then the denominator may be zero if and only if $x$ is a common zero of $P_n$ and $P_{n-1}$. But this is impossible since $P_n$ and $P_{n-1}$ have no common zeros~\cite{1}.

Thus from the limit expression~(\ref{limit_R_n}) it follows readily that the absolutely continuous parts of the spectrum of $A_n$ and $J_n$ coincide. Prove that under the conditions of the theorem the singular part of the spectrum of $A_n$ is absent. For this purpose, return to the formula~(\ref{R_n}). The numerator and the denominator in this formula are analytic functions on the set
$\mathbb C\setminus\sigma(J_n)$. Therefore the point $x\in\mathbb R$ belongs to the singular part of the spectrum of $A_n$ for $|x-a_n|\le2b_n$ if and only if
\begin{equation}
\label{3.3}
\begin{array}{l}
\lim\limits_{\epsilon\to+0} (P_n(x+i\epsilon)+P_{n-1}(x+i\epsilon)\,b_{n-1}\,K_n(x+i\epsilon))=\\
=P_n(x)+P_{n-1}(x)\,b_{n-1}D_n(x)+i\,P_{n-1}(x)\,b_{n-1}B_n(x)=0
\end{array}
\end{equation}
and for $|x-a_n|>2b_n$ if and only if $x$ is a zero of the denominator:
\begin{equation}
\label{3.4}
P_n(x)+P_{n-1}(x)\,b_{n-1}K_n(x)=P_n(x)+P_{n-1}(x)\,b_{n-1}D_n(x)=0
\end{equation}
But as it was shown before, the equality~(\ref{3.3}) is impossible for $|x-a_n|<2b_n$. It remains to consider the case $|x-a_n|\ge2b_n$. In this case the conditions~(\ref{3.3}) and~(\ref{3.4}) coincide since $B_n(x)=0$.

Thus, let $|x-a_n|\ge2b_n$. From the definition (\ref{limit_K_n2}) of $D_n(x)$ it follows that $\forall x\in\mathbb R$
$$
|D_n(x)|\le\frac{1}{b_n}
$$
Using this estimate and the assumed monotonicity of the sequence $b_n$ one obtains
\begin{equation}
\label{3.5}
\begin{array}{l}
|P_n(x)+P_{n-1}(x)\,b_{n-1}D_n(x)|\ge|P_n(x)|-|P_{n-1}(x)|b_{n-1}|D_n(x)|\ge\\
\ge|P_n(x)|-|P_{n-1}(x)|\,\frac{\displaystyle b_{n-1}}{\displaystyle b_n}\ge |P_n(x)|-|P_{n-1}(x)|
\end{array}
\end{equation}
On the other hand, from the recurrence relations~(\ref{recurrent_relations}) for the polynomials $P_n$ one has
$$
|P_{n}(x)|=\left|\frac{(x-a_{n-1})P_{n-1}(x)-b_{n-2}P_{n-2}(x)}{b_{n-1}}\right|\ge\frac{|x-a_{n-1}|}{b_{n-1}}|P_{n-1}(x)|-
\frac{b_{n-2}}{b_{n-1}}|P_{n-2}|
$$
If $|x-a_n|\ge2b_n$, then by virtue of the theorem condition $|a_n-a_{n-1}|\le2(b_n-b_{n-1})$, the inequality $|x-a_{n-1}|\ge2b_{n-1}$ also holds. Therefore from the last inequality one obtains
$$
|P_{n}(x)|\ge2|P_{n-1}(x)|-|P_{n-2}(x)|
$$
or
$$
|P_{n}(x)|-|P_{n-1}(x)|\ge|P_{n-1}(x)|-|P_{n-2}(x)|
$$
Continuing this process by induction one obtains
$$
|P_{n}(x)|-|P_{n-1}(x)|\ge|P_{1}(x)|-|P_{0}(x)|\ge2-1=1
$$
Finally, substituting this inequality into~(\ref{3.5}) one arrives at the following estimate of the denominator in the formula~(\ref{R_n}) in the case
$|x-a_n|\ge2b_n$
$$
|P_n(x)+P_{n-1}(x)\,b_{n-1}D_n(x)|\ge1\,,
$$
whence it follows the absence in this area the points of a singular spectrum of $A_n$.

It remains to prove the formula~(\ref{f_n(x)_3.6}). Using~(\ref{f_J_n}) and (\ref{limit_R_n}) one has
$$
f_n(x)=\frac{1}{\pi}\lim\limits_{\epsilon\to+0}\mbox{Im}\,R_n(x+i\epsilon)=\frac{\frac{\displaystyle 1}{\displaystyle \pi}\,B_n(x)}
{(P_n(x)+b_{n-1}\,P_{n-1}(x)\,D_n(x))^2+(b_{n-1}\,P_{n-1}(x)\,B_n(x))^2}=
$$
$$
=\frac{f_{J_n}(x)}{(P_n(x)+b_{n-1}\,P_{n-1}(x)\,D_n(x))^2+(b_{n-1}\,P_{n-1}(x)\,B_n(x))^2}
$$
Transform the denominator. When $|x-a_n|\le2b_n$, using~(\ref{limit_K_n2}), (\ref{limit_K_n3}), (\ref{recurrent_relations}), one has
$$
(P_n(x)+b_{n-1}\,P_{n-1}(x)\,D_n(x))^2+(b_{n-1}\,P_{n-1}(x)\,B_n(x))^2=
$$
$$
=\left(P_n+\frac{a_n-x}{2b_n^2}b_{n-1}P_{n-1}\right)^2
+P_{n-1}^2b_{n-1}^2\frac{4b_n^2-(a_n-x)^2}{4b_n^4}=
$$
$$
=P_n^2+\frac{a_n-x}{b_n^2}b_{n-1}P_{n-1}P_n+\frac{b_{n-1}^2}{b_n^2}P_{n-1}^2=P_n^2+P_{n-1}\frac{b_{n-1}}{b_n^2}((a_n-x)P_n+
b_{n-1}P_{n-1})=
$$
$$
=P_n^2-\frac{b_{n-1}}{b_n}P_{n-1}P_{n+1}
$$
The theorem is proved.\\
-------------------------------------------------------------------------------

{\bf Remark.}\, Additional condition of the theorem means some subordination of $a_n$ to $b_n$. Indeed, for $k\in{\mathbb N}$, one has
$$
|a_k|-|a_{k-1}|\le|a_k-a_{k-1}|\le2(b_k-b_{k-1})
$$
Adding these inequalities for $k=1,2,\ldots,n$ one obtains
$$
|a_n|\le|a_0|+2(b_n-b_0)
$$
-------------------------------------------------------------------------------

Note also that the expression 
$$
P_n^2(x)-\frac{b_{n-1}}{b_n}P_{n-1}(x)P_{n+1}(x)
$$
appears probably first time in works~\cite{20}. In~\cite{21} it is called "Turan determinant".

\vspace{0.8cm}

{\bf Definition.}\,
Let us call the system of intervals $I_n$ {\it centered} on some interval $[a,b\,]$ if this interval is contained inside the interval $I_n$ for all sufficiently large values of~$n$.\\

-------------------------------------------------------------------------------

Now, when the absolute continuity conditions for the spectrum of $A_n$ are stated, we can apply the theorem~(\ref{criterion_of_absolutely_continuity}). Combining the results of the theorems~(\ref{criterion_of_absolutely_continuity}) and~(\ref{spectr_A_n}) and using the Definition, one obtains

\begin{theorem}
\label{Jacobi_criterion_of_absolutely_continuity}
Assume that the operator $A$, which is represented by~(\ref{Jacobi_Matrix}), is self-adjoint and the system of intervals $I_n$ is centered on some interval $[a,b\,]$ of the operator $A$ spectrum. Then the spectrum of the operator $A$ is absolutely continuous on $[a,b\,]$ if on this interval the functions
$$
f_n(x)=\frac{f_{J_n}(x)}{P^2_n(x)-\frac{\displaystyle b_{n-1}}{\displaystyle b_n}P_{n-1}(x)P_{n+1}(x)}
$$
have equi-absolutely continuous integrals for all $n$ sufficiently large.
\end{theorem}

One can give a simple sufficient condition of equi-absolute continuity of functions $\sigma_n(\lambda)$ in the case when
$f_n(x)\in L_\infty[a,b\,]$ for all sufficiently large $n$, that is, when $\|f_n(x)\|_{L_\infty}=vrai\,sup\,\{f_n(x)\}<+\infty$. Let the norms of functions $f_n(x)$ be uniformly bounded in the space $L_\infty[a,b\,]$
$$
\|f_n(x)\|_{L_\infty}<C\,,
$$
where the constant $C$ is independent of $n$. Then for any measurable set $e\subseteq[a,b\,]$ one has
$$
\int\limits_{e}f_n(t)\,dt<C\,m(e)
$$
Because the constant $C$ is independent of $n$, for any $\epsilon>0$ and $m(e)<\epsilon/C$ it follows
$$
\int\limits_{e}f_n(t)\,dt<\epsilon
$$
This means that the functions $f_n(x)$ have equi-absolutely continuous integrals.

In the case of the operators $A_n$, if the system of the intervals $I_n$ is centered on the interval $[a,b\,]$, then for all sufficiently large values of $n$ the functions $f_n(x)$ are continuous on $[a,b\,]$ and $vrai\,sup\,\{f_n(x)\}=\max\limits_{x\in[a,b]}\{f_n(x)\}$. Thus we have the following sufficient condition for the absolute continuity of the spectrum of operator $A$:
\begin{theorem}
\label{absolute_continuity_A_1}
Assume that the operator $A$, which is represented by~(\ref{Jacobi_Matrix}), is self-adjoint and the system of intervals $I_n$ is centered on some interval $[a,b\,]$ of the operator $A$ spectrum. Let for all sufficiently large $n$
\begin{equation}
\label{condition_for_f_n(x)}
\max\limits_{x\in[a,b]}\{f_n(x)\}<C\,,
\end{equation}
where the constant $C>0$ is independent of $n$. Then the spectrum of $A$ is absolutely continuous on $[a,b\,]$ and the sequence  $\sigma_n(x)$ uniformly converges to $\sigma(x)$.
\end{theorem}

An expression~(\ref{f_n(x)_3.6}) for the spectral wights $f_n(x)$ can be transformed to another form which is sometimes more convenient for analysis. Use for that the following assertion~( see also \cite{10}).

\begin{lemma}
\label{f_n(x)_3.7}
Let $b_{-1}\equiv0$, $P_{-1}(x)\equiv0$. Then for $n=0,1,2,\ldots$ the identity
$$
P^2_n(x)-\frac{b_{n-1}}{b_n}P_{n-1}(x)P_{n+1}(x)=\frac{1}{b_n^2}\sum\limits_{k=0}^n\left[(b_k^2-b_{k-1}^2)P_k^2(x)+
b_{k-1}(a_k-a_{k-1})P_{k-1}(x)P_k(x)\right]
$$
\end{lemma}
is valid.

{\bf Proof.}\, Consider two successive recurrence relations~(\ref{recurrent_relations}) for polynomials $P_n(x)$
$$
b_{k-1}P_{k}=(x-a_{k-1})P_{k-1}-b_{k-2}P_{k-2}
$$
$$
b_kP_{k+1}=(x-a_{k})P_{k}-b_{k-1}P_{k-1}
$$
Multiplying the first of them by $P_k$ and the second by $P_{k-1}$ and subtracting, one obtains
$$
(a_k-a_{k-1})P_{k-1}P_k=b_{k-1}(P_k^2-P_{k-1}^2)+b_{k-2}P_{k-2}P_k-b_kP_{k-1}P_{k+1}
$$
Multiplying this equality by $b_{k-1}$ and adding for $k$ from 1 up to $n$, one finds
$$
\sum\limits_{k=1}^nb_{k-1}(a_k-a_{k-1})P_{k-1}P_k=\sum\limits_{k=1}^nb_{k-1}^2(P_k^2-P_{k-1}^2)-b_{n-1}b_nP_{n-1}P_{n+1}
$$
or
$$
\sum\limits_{k=0}^{n-1}(b_k^2-b_{k-1}^2)P_k^2+\sum\limits_{k=1}^{n}b_{k-1}(a_k-a_{k-1})P_{k-1}P_k=b_{n-1}^2P_n^2-
b_{n-1}b_nP_{n-1}P_{n+1}
$$
or
$$
\sum\limits_{k=0}^{n}(b_k^2-b_{k-1}^2)P_k^2+\sum\limits_{k=0}^{n}b_{k-1}(a_k-a_{k-1})P_{k-1}P_k=b_{n}^2P_n^2-
b_{n-1}b_nP_{n-1}P_{n+1}
$$
Dividing this equality by $b_n^2$, one obtains the required statement.\\
-------------------------------------------------------------------------------

Using this lemma one obtains

\begin{theorem}
\label{absolute_continuity_A}
Assume that the operator $A$, which is represented by~(\ref{Jacobi_Matrix}), is self-adjoint and the system of intervals $I_n$ is centered on some interval $[a,b\,]$ of the operator $A$ spectrum. If for any $n\ge N$ the estimate
\begin{equation}
\label{main_estimation}
\sum\limits_{k=0}^n\left[(b_k^2-b_{k-1}^2)P_k^2(x)+
b_{k-1}(a_k-a_{k-1})P_{k-1}(x)P_k(x)\right]>Cb_n\,,
\end{equation}
holds, where the constant $C>0$ depends neither on $n$ no on $x\in[a,b\,]$, then the spectrum of $A$ is absolutely continuous on $[a,b\,]$ and the sequence $\sigma_n(x)$ uniformly converges on this interval to the function $\sigma(x)$.
\end{theorem}

{\bf Proof.}\, By virtue of the Theorem~(\ref{absolute_continuity_A_1}), for the proof it is remained to show that~(\ref{condition_for_f_n(x)}) follows from~(\ref{main_estimation}). From~(\ref{f_J_n}) and (\ref{f_n(x)_3.6}) one has
$$
f_n(x)\le\frac{1/(\pi b_n)}{P^2_n(x)-\frac{\displaystyle b_{n-1}}{\displaystyle b_n}P_{n-1}(x)P_{n+1}(x)}\,,
$$
for any $x$ and $n$, or by the Lemma~(\ref{f_n(x)_3.7})
$$
f_n(x)\le\frac{b_n}{\pi\sum\limits_{k=0}^n\left[(b_k^2-b_{k-1}^2)P_k^2(x)+
b_{k-1}(a_k-a_{k-1})P_{k-1}(x)P_k(x)\right]}
$$
Hence, if~(\ref{main_estimation}) holds, then $f_n(x)<(\pi C)^{-1}$. Theorem is proved.\\
-------------------------------------------------------------------------------

\begin{corollary}
\label{corollary_a_n}
If $a_n=0$, $n=0,1,\ldots$, then the condition~(\ref{main_estimation}) takes the form
$$
\sum\limits_{k=0}^n(b_k^2-b_{k-1}^2)P_k^2(x)>Cb_n\,,
$$
\end{corollary}

Let us give one more variant of sufficient condition for absolute continuity based on one more representation of $f_n(x)$ which is used below.

Using the recurrence relations for polynomials $P_n$, one can present the denominator of functions $f_n(x)$ in the form
$$
P^2_n(x)-\frac{\displaystyle b_{n-1}}{\displaystyle b_n}P_{n-1}(x)P_{n+1}(x)=P_{n+1}^2(x)+P_{n}^2(x)-\frac{x-a_n}{b_n}P_{n+1}(x)P_{n}(x)
$$
Hence for $|x-a_n|\le b_n$ the estimate
$$
P^2_n(x)-\frac{\displaystyle b_{n-1}}{\displaystyle b_n}P_{n-1}(x)P_{n+1}(x)\ge \frac{1}{2}\left(P_{n+1}^2(x)+P_{n}^2(x)\right)
$$
holds. Therefore
$$
f_n(x)\le\frac{2/(\pi b_n)}{P_{n+1}^2(x)+P_{n}^2(x)}\,,
$$
whence one obtains

\begin{theorem}
\label{absolute_continuity_A_2}
Suppose that the operator $A$ represented by Jacobi matrix~(\ref{Jacobi_Matrix}) is self-adjoint and the domain $|x|\le K$ belongs to its spectrum. If for any $x$ from the domain $|x|\le K$ and any $n\ge N$ the estimates $|x-a_n|\le b_n$ and
$$
P_{n+1}^2(x)+P_{n}^2(x)>\frac{C}{b_n}
$$
hold, where the constant $C>0$ depends neither on $x$ no on $n$, then the spectrum of the operator $A$ is absolutely continuous in the domain $|x|\le K$ and the sequence $\sigma_n(x)$ uniformly converges in this domain to the function $\sigma(x)$.
\end{theorem}
-------------------------------------------------------------------------------

{\bf Examples.}\, Let us consider some examples. Consider first the case $a_n=0$. In this case
$I_n=[-2b_n,2b_n]$. Let the sequence $b_n$ be non-decreasing. Two cases are possible: $b_n$ tends either to a finite limit or to the infinity. Let's start with the first case as the most simple.

Let $b_n\to b<+\infty$. In this case the operator $A$ is bounded. The difference $A-J$, where
$$
J=\left(\begin{array}{ccccc} 0&b&0&0&\ldots\\
                      b&0&b&0&\ldots\\
                      0&b&0&b&\ldots\\
                      0&0&b&0&\ldots\\
                      \vdots&\vdots&\vdots&\vdots&\ddots
  \end{array}\right)\,,
$$
is compact and Weyl's theorem gives at once $\sigma_{ess}(A)=\sigma_{ess}(J)=[-2b,2b\,]$.

The functions $f_n(x)$ have the form
$$
f_n(x)=\frac{f_{J_n}(x)}{P^2_n(x)-\frac{\displaystyle b_{n-1}}{\displaystyle b_n}P_{n-1}(x)P_{n+1}(x)}=
\frac{b_n^2\,f_{J_n}(x)}{\sum\limits_{k=0}^n(b_k^2-b_{k-1}^2)P_k^2(x)}\,,
$$
where
$$
f_{J_n}(x)=\left[\begin{array}{cr}\frac{\displaystyle\sqrt{4(b_n)^2-x^2}}{\displaystyle 2\pi\,(b_n)^2}\,,& |x|\le2b_n\\
                             0\,,& |x|>2b_n
             \end{array}\right.\,,
$$
A simple estimate gives
$$
f_n(x)\le\frac{b}{\pi b_0^2}
$$
for all $n$ and $x$. Hence the functions $f_n(x)$ are uniformly bounded and the Theorem~(\ref{absolute_continuity_A_1}) gives at once the absolute continuity of the spectrum of $A$ on any interval inside $[-2b,2b\,]$ . Thus the absolutely continuous spectrum
is concentrated on the interval $[-2b,2b\,]$.

Let's consider the case $b_n\to+\infty$. The system of intervals $I_n$ is centered now on any interval of a real axis.
Take as an example Jacobi matrix with elements $a_n=0$ and $b_n=\sqrt{\frac{n+1}{2}}$, $n=0,1,\ldots$. This matrix is connected with Hermite polynomials. The corresponding operator is self-adjoint by virtue of Carleman's condition~(\ref{Charleman's_condtion}). In this case $b_k^2-b_{k-1}^2=1/2$ for any $k$ and the function $f_n(x)$ can be presented in the form
$$
f_n(x)=\frac{(n+1)\,f_{J_n}(x)}{\sum\limits_{k=0}^nP_k^2(x)}
$$
The orthogonal polynomials of 1st type coincide with Hermit polynomials with accuracy to a constant~\cite{13}
$$
P_n(x)=\frac{1}{\sqrt{n!\,2^n}}\,H_n(x)\,,\quad H_n(x)=(-1)^n\,e^{x^2}\frac{d^{\,n}}{dx^n}\,e^{-x^2}
$$
The functions $f_n(x)$ are continuous and have a maximum at some points inside the interval $|x|\le 2b_n$. Show that this maximum is attained by all of this functions at the point $x=0$, that is,
$$
f_n(x)\le f_n(0)
$$
Consider separately the numerator and the denominator in the formula for $f_n(x)$. The numerator, that is the function $f_{J_n}(x)$, has a maximum at the point $x=0$. Show that the denominator is minimal at the same point. Consider for this the derivative of the expression $P_k^2(x)+P_{k-1}^2(x)$:
$$
\left(P_k^2(x)+P_{k-1}^2(x)\right)'=2P_k(x)\,P_k\,'(x)+2P_{k-1}(x)\,P_{k-1}'(x)
$$
Using the expression for the derivative of Hermite polynomials~\cite{8}
$$
H_k'(x)=2kH_{k-1}(x)\,,
$$
one obtains
$$
P_k'(x)=\sqrt{2k}P_{k-1}(x)\,,
$$
and hence
$$
\left(P_k^2(x)+P_{k-1}^2(x)\right)'=2\sqrt{2}\,P_{k-1}\left(\sqrt{k}\,P_k+\sqrt{k-1}\,P_{k-2}\right)
$$
From the recurrence relation for $P_k$
$$
\sqrt{\frac{k}{2}}\,P_k=x\,P_{k-1}-\sqrt{\frac{k-1}{2}}\,P_{k-2}\,,
$$
one obtains
$$
\left(P_k^2(x)+P_{k-1}^2(x)\right)'=4x\,P_{k-1}^2(x)
$$
Using this relation, we can present the derivative of the denominator in the expression for $f_n(x)$ for even and odd values of $n\ge 1$ in the form
\begin{equation}
\label{derivatives}
\left(\sum\limits_{k=0}^{2m-1}P_k^2(x)\right)'=4x\left(P_0^2(x)+P_2^2(x)+\ldots+P_{2m-2}^2(x)\right)\,,\quad m=1,2,\ldots
\end{equation}
$$
\left(\sum\limits_{k=0}^{2m}P_k^2(x)\right)'=4x\left(P_1^2(x)+P_3^2(x)+\ldots+P_{2m-1}^2(x)\right)\,,\quad m=1,2,\ldots
$$
From this formulaes it follows readily that for $x>0$ the expression $\sum\limits_{k=0}^{n}P_k^2(x)$ increases and for $x<0$ decreases. Thus the minimal value of $\sum\limits_{k=0}^{n}P_k^2(x)$ is attained at the point $x=0$ and hence
$$
f_n(x)\le f_n(0).
$$
It follow that in order that the sequence of functions $f_n(x)$ be uniformly bounded it is sufficient that the number sequence $f_n(0)$ be bounded above. It is not difficult to do not only in this specific example but in a much more general case.

Indeed for $a_n=0$ and $x=0$ the recurrence relations~(\ref{recurrent_relations}) take the form
$$
b_n\,P_{n+1}=-b_{n-1}\,P_{n-1}\,,
$$
whence, using initial conditions $P_0\equiv1$, $P_1(0)=0$, one obtains
$$
P_1(0)=P_3(0)=\ldots=P_{2k-1}(0)=\ldots=0
$$
$$
P_{2k}(0)=(-1)^k\,\frac{b_{2k-2}b_{2k-4}\cdot\ldots\cdot b_2\,b_0}{b_{2k-1}b_{2k-3}\cdot\ldots\cdot b_3\,b_1}
$$
Substituting this relations into the expression for $f_n(0)$
$$
f_n(0)=\frac{f_{J_n}(0)}{P^2_n(0)-\frac{\displaystyle b_{n-1}}{\displaystyle b_n}P_{n-1}(0)P_{n+1}(0)}\,,
$$
one finds for even and odd values of $n$
\begin{equation}
\label{f_n(0)}
f_{2k-1}(0)=\frac{1}{\pi}\,\frac{1}{b_{2k-1}}\,\frac{b_{2k-1}^2b_{2k-3}^2\cdot\ldots\cdot b_3^2\,b_1^2}{b_{2k-2}^2b_{2k-4}^2\cdot\ldots
\cdot b_2^2\,b_0^2}
\end{equation}
$$
f_{2k}(0)=\frac{1}{\pi}\,\frac{1}{b_{2k}}\,\frac{b_{2k-1}^2b_{2k-3}^2\cdot\ldots\cdot b_3^2\,b_1^2}{b_{2k-2}^2b_{2k-4}^2\cdot\ldots\cdot
b_2^2\,b_0^2}
$$
Since $b_{2k}\ge b_{2k-1}$,
$$
f_{2k}(0)\le f_{2k-1}(0)\,,
$$
and the boundedness of $f_n(0)$ reduces to the boundedness of their odd subsequence.

From
$$
\frac{f_{2k+1}(0)}{f_{2k-1}(0)}=\frac{b_{2k-1}b_{2k+1}}{b_{2k}^2}
$$
it follows that if
\begin{equation}
\label{monoton}
b_{2k-1}b_{2k+1}\le b_{2k}^2\,,\quad k=1,2,\ldots\,,
\end{equation}
then the sequence $f_{2k-1}(0)$ is non-increasing. A wide enough class of sequences $b_n$ satisfies this condition. In particular, all power sequences $b_n=C\,(n+1)^a$, $a\in(0;1]$, $C>0$ belong to this class.
The condition $a\in(0,1]$ does not restrict the validity of~(\ref{monoton}) but for $a>1$ the operator $A$ corresponding to Jacobi matrix become not self-adjoint by Berezanskii's theorem~\cite{1} since in this case the series $\sum\frac{1}{b_n}$ diverges. But any self-adjoint extension of such an operator has a purely discrete spectrum~\cite{13}.

By the same way one can show that if the condition~(\ref{monoton}) is fulfilled for all $n$, that is $b_{n-1}b_{n+1}\le b_{n}^2$,
then together with non-increase of $f_{2k-1}(0)$, the even subsequence $f_{2k}(0)$ is non-decreasing.

Returning to the above example $a_n=0$, $b_n=\sqrt{\frac{n+1}{2}}$, by the Theorem~(\ref{absolute_continuity_A_1}) one concludes that the spectrum of the corresponding operator $A$ may be absolutely continuous only on any interval of a real axis and hence everywhere since $\max\limits_x\{f_n(x)\}=f_n(0)\le f_{2k-1}(0)\le f_1(0)$. Since the spectrum of self-adjoint operator is not empty the spectrum of $A$ is purely absolutely continuous. It is well known~\cite{13} that in this case the distribution function of the operator $A$ is defined by
$$
\sigma(\lambda)=\frac{1}{\sqrt\pi}\int\limits_{-\infty}^{\lambda}e^{-x^2}\,dx\,,
$$
and the absolutely continuous spectrum covers the whole real axis.

If we take the sequence $b_n$ of the form $b_n=\left(\frac{n+1}{2}\right)^a$ where $a\in(\frac{1}{2};1]$, then the direct calculations show that the behavior of the expression $\sum\limits_{k=0}^n(b_k^2-b_{k-1}^2)P_k^2(x)$ for any $n$ is the same as in the case of Hermite polynomials at $a=1/2$. That is, it has a minimum at $x=0$, increasing for $x>0$ and decreasing for $x<0$. All that is said above leads to the conclusion that in this case the operator $A$ has purely absolutely continuous spectrum too.

In the case $a\in (0;\frac{1}{2})$ the situation changes and the function $\sum\limits_{k=0}^n(b_k^2-b_{k-1}^2)P_k^2(x)$ can attain a minimum at $x\ne 0$.

The sequence $f_n(0)$ is not always bounded. Indeed, saving a monotonicity of a sequence $b_n$, let
$b_{2k-1}=b_{2k}\,,\:k=1,2,\ldots$. From (\ref{f_n(0)}) one obtains
$$
f_{2k-1}(0)=f_{2k}(0)=\frac{b_{2k-1}}{\pi b_0^2}
$$
Thus if $b_{n}\to\infty$, then $f_n(0)\to\infty$. If this singularity is strong enough, then in the point $x=0$ an eigenvalue of the operator $A$ appears. Indeed in this case
$$
P^2_{2k}(0)=\frac{b_0^2}{(b_{2k-1})^2}\,,
$$
and if the series
$$
\sum_k\frac{1}{(b_{2k-1})^2}
$$
converges, then so is $\sum\limits_n P^2_n(0)$. It means that the point $x=0$ is an eigenvalue.

For the spectral analysis it is of great interest to find more or less general conditions for the sequences $a_n$ and $b_n$ which guarantee the absolute continuity of the spectrum of $A$. State below one of such conditions (see also~\cite{19,22,23})

\begin{theorem}
\label{absolute_continuity_A_trough_a_n_b_n}
Let $a_n$ and $b_n$ satisfy the conditions
\begin{enumerate}
\item $\lim b_n=+\infty\,,\qquad \sum\frac{\displaystyle 1}{\displaystyle b_n}=+\infty$

\item $\displaystyle 0<C_1<\frac{b_{n+1}}{b_n}<C_2<+\infty$

\item $\{\frac{\displaystyle a_n}{\displaystyle b_n}\}\in l_2$, if $\lim|a_n|=+\infty$ or\\
$\{\frac{\displaystyle 1}{\displaystyle b_n}\}\in l_2$, if $|a_n|$ is bounded

\item $\displaystyle\left\{\frac{b_{n-1}}{b_n}-\frac{b_{n-2}}{b_{n-1}}\right\}
\in l_1\,,\quad
\left\{\frac{b_n-b_{n-2}}{b_{n-1}b_n}\right\}\in l_1\,,\quad
\left\{\frac{b_na_{n-1}-a_nb_{n-2}}{b_{n-1}b_n}\right\}\in l_1$

\item $\displaystyle\underline{\lim}\,\frac{b_{n-1}^2b_{n-3}^2b_{n-5}^2\ldots}{b_nb_{n-2}^2b_{n-4}^2\ldots}>0$
\end{enumerate}
Then the spectrum of the operator $A$ is purely absolutely continuous and the sequence $\sigma_n(x)$ uniformly converges to the function $\sigma(x)$ on any finite interval.
\end{theorem}

{\bf Proof.}\, First, from the conditions (1) of the theorem it follows that the operator $A$ is self-adjoint. Let $K$ be an arbitrarily large positive number. By virtue of the conditions (1) and (3) of a theorem for any $x$ such that $|x|\le K$, the inequality $|x-a_n|\le b_n$ holds for all sufficiently large $n$. To apply the Theorem~(\ref{absolute_continuity_A_2}) it remains to prove uniform estimate
$$
P_{n+1}^2(x)+P_{n}^2(x)>\frac{C}{b_n}
$$
It is convenient to use for that the apparatus of $(2\times2)$ matrices. Three term recurrence relations for the polynomials can be written down in the matrix form
$$
\vec{u}_{n+1}=B_n\vec{u}_n\,,\quad n=1,2,\ldots\,,
$$
where
$$
\vec{u}_n=\left(\begin{array}{c} P_{n-1} \\
                                 P_n
                \end{array}\right)\,,\qquad B_n=\left(\begin{array}{cc} 0 & 1\\
                                                                       -\frac{b_{n-1}}{b_n} & \frac{x-a_n}{b_n}
                \end{array}\right)
$$
From this it follows that
$$
\vec{u}_{n+1}=B_nB_{n-1}\ldots B_1\vec{u}_1
$$
and
$$
P_{n+1}^2+P_n^2=\|\vec{u}_{n+1}\|^2
$$
To obtain the required estimate of this expression from below, we use the idea of the work~\cite{15}. One has
$$
B_nB_{n-1}=\left(\begin{array}{cc} -\frac{b_{n-2}}{b_{n-1}} & \frac{x-a_{n-1}}{b_{n-1}}\\
                                   -\frac{b_{n-2}}{b_{n-1}}\frac{x-a_n}{b_n} & -\frac{b_{n-1}}{b_{n}}+\frac{(x-a_n)(x-a_{n-1})}
{b_nb_{n-1}}
                \end{array}\right)
$$
By virtue of the conditions (2) and (3) of the theorem one obtains
$$
B_nB_{n-1}=\left(\begin{array}{cc} -\frac{b_{n-2}}{b_{n-1}} & \frac{x-a_{n-1}}{b_{n-1}}\\
                                   -\frac{b_{n-2}}{b_{n-1}}\frac{x-a_n}{b_n} & -\frac{b_{n-1}}{b_{n}}
                \end{array}\right)+R_n\,,
$$
where $\{\|R_n\|\}\in l_1$ and the series $\sum\|R_n\|$ converges uniformly in the domain $|x|\le K$. Next one has
$$
B_nB_{n-1}=\left(\begin{array}{cc} -\frac{b_{n-2}}{b_{n-1}} & 0\\
                                   0 & -\frac{b_{n-1}}{b_{n}}
                \end{array}\right)+\left(\begin{array}{cc} 0 & \frac{x-a_{n-1}}{b_{n-1}}\\
                                   -\frac{b_{n-2}}{b_{n-1}}\frac{x-a_n}{b_n} & 0
                \end{array}\right)+R_n\,,
$$
and by virtue of the condition (4)
$$
B_nB_{n-1}=-\frac{b_{n-1}}{b_{n}}I+\rho_n P+R_n'\,,\quad P=\left(\begin{array}{cc} 0& 1\\
                                                                                  -1 & 0
                \end{array}\right)\,,\quad \rho_n=\frac{x-a_n}{b_n}
$$
Here also as above $\{\|R_n'\|\}\in l_1$ uniformly in the domain $|x|\le K$. One has next
$$
B_nB_{n-1}=-\frac{b_{n-1}}{b_{n}}\left(I-\frac{b_{n}}{b_{n-1}}\rho_n P+R_n''\right)=
-\frac{b_{n-1}}{b_{n}}\left(e^{-\frac{b_{n}}{b_{n-1}}\rho_n P}+R_n'''\right)
$$
$$
=-\frac{b_{n-1}}{b_{n}}e^{-\frac{b_{n}}{b_{n-1}}\rho_n P}\left(I+S_n\right)
$$
Here the matrix $e^{-\frac{b_{n}}{b_{n-1}}\rho_n P}$ is orthogonal and as above $\{\|S_n\|\}\in l_1$ and the corresponding series converges uniformly in the domain $|x|\le K$. Using this for all sufficiently large $n$, such that $\|S_n\|<1$, one obtains
$$
P_{n+1}^2+P_n^2\ge\frac{b_{n-1}^2b_{n-3}^2\ldots b_{N-1}^2}{b_n^2b_{n-2}^2\ldots b_N^2}\prod\limits_{k=0}^{(n-N)/2}(1-\|S_{n-2k}\|)^2\,
\|\vec{u}_N\|^2
$$
From this by virtue of the convergence of an infinite product and the condition (5) of the theorem one obtains uniformly in the domain $|x|\le K$
$$
P_{n+1}^2+P_n^2\ge\frac{C}{b_n}\,,
$$
By Theorem~(\ref{absolute_continuity_A_2}) the spectrum of the operator $A$ may be absolutely continuous only in the domain $|x|\le K$. Since the number $K$ is arbitrarily large the same is valid for any interval. Since the spectrum of self-adjoint operator is not empty the spectrum of $A$ is purely absolutely continuous. Theorem is proved.\\
-------------------------------------------------------------------------------

{\bf Remark.}\, The first condition of the item 4 of the theorem and the condition of the item 5 are fulfilled if $b_n^2\ge b_{n-1}b_{n+1}$. The last as it was noticed already above is valid for all sequences of the form $b_n=(n+1)^a$ for $a\in(0;1]$.\\
-------------------------------------------------------------------------------

\section{Convergence of the derivatives of distribution functions.}

In the previous chapter it was shown that under some conditions the sequence of the absolutely continuous on the interval $[a,b\,]$ functions $\sigma_n(\lambda)$, associated with the operators $A_n$,
$$
\sigma_n(\lambda)=\sigma_n(a)+\int\limits_{a}^\lambda f_n(x)\,dx
$$
converges at every point $\lambda\in[a,b\,]$ to the absolutely continuous function $\sigma(\lambda)$
$$
\sigma(\lambda)=\sigma(a)+\int\limits_{a}^\lambda f(x)\,dx\,,
$$
associated with the operator $A$. Moreover the convergence is uniform on $[a,b\,]$.

The following questions arise: what can be said about the convergence of a sequence of the derivatives $f_n(x)$? Does this sequence converge to $f(x)$? And if it converges, then in what sense?

Let the system of intervals $I_n$ be centered on some interval $[a,b\,]$ of the spectrum of $A$ which does not contain the eigenvalues. Note that we do not assume that the spectrum of $A$ is absolutely continuous on $[a,b\,]$. Consider for sufficiently large $n$ a sequence of absolutely continuous functions $\sigma_n(\lambda)$ on this interval. From the results of the previous chapter (Lemma~\ref{lemma2}) it follows that $\sigma_n(\lambda)\to\sigma(\lambda)$ for all $\lambda\in[a,b\,]$. The question about the existence of the limit of sequence $f_n(x)$ and its equality to $\sigma'(x)$ is equivalent to the question about the possibility of changing the order of two limits. Indeed by virtue of a continuity of the functions $f_n(x)$ we can write down
\begin{equation}
\label{two_limits_1}
\lim\limits_{n\to\infty}f_n(x)=\lim\limits_{n\to\infty}\sigma_n'(x)=\lim\limits_{n\to\infty}\,\lim\limits_{x'\to x}
\frac{\sigma_n(x)-\sigma_n(x')}{x-x'}
\end{equation}
On the other hand, assuming the existence of a derivative $\sigma'(x)$, one obtains
\begin{equation}
\label{two_limits_2}
\sigma'(x)=\lim\limits_{x'\to x}\frac{\sigma(x)-\sigma(x')}{x-x'}=\lim\limits_{x'\to x}\,\lim\limits_{n\to\infty}
\frac{\sigma_n(x)-\sigma_n(x')}{x-x'}
\end{equation}
This expression differs from~(\ref{two_limits_1}) by changing the order of limits only. We know that internal limits in~(\ref{two_limits_1}) and ~(\ref{two_limits_2}) exist. In this case, as known, for the existence and equality finite external limits it is sufficient that one of the internal limits be uniform~\cite{14}. Consider internal limit in~(\ref{two_limits_1}). Since the functions $f_n(x)$ are continuous by mean value theorem one has
$$
\lim\limits_{x'\to x}\frac{\sigma_n(x)-\sigma_n(x')}{x-x'}=\lim\limits_{x'\to x}\frac{1}{x-x'}\,\int\limits_{x'}^x f_n(t)\,dt=
\lim\limits_{x'\to x}f_n(\xi_n)\,,
$$
where $\xi_n$ lies between $x$ and $x'$ and uniformly tends to $x$ as $x'\to x$. From this, one sees, that in order that this limit be uniform it is sufficient to require that convergence $f_n(x')$ to $f_n(x)$ be uniform in $n$. If it is fulfilled for all $x$ in $[a,b\,]$, then the sequence $f_n(x)$ is equicontinuous on this interval. But as known~\cite{12}, from the existence of finite derivative $\sigma'(x)$ everywhere in $[a,b\,]$ it follows the absolute continuity of $\sigma(x)$ on this interval. Hence we obtain
\begin{theorem}
\label{absolute_continuity_through_uniform_continuity}
Let the system of intervals $I_n$ be centered on some interval $[a,b\,]$ of the operator $A$ spectrum which does not contain the eigenvalues. Suppose that for all sufficiently large $n$ the sequence of functions $f_n(x)$ is equicontinuous on this interval. Then the spectrum of $A$ is absolutely continuous on $[a,b\,]$ and $\forall x\in[a,b\,]$
$$
\lim\limits_{n\to\infty} f_n(x)=\sigma'(x)
$$
\end{theorem}
-------------------------------------------------------------------------------

Combining this theorem with the theorem~(\ref{absolute_continuity_A_1}) and using Arzela-Ascoli theorem one obtains
\begin{theorem}
\label{absolute_continuity_and_convergence}
Assume that the operator $A$ represented by Jacobi matrix~(\ref{Jacobi_Matrix}) is self-adjoint and the system of intervals $I_n$ is centered on some interval $[a,b\,]$ of the operator $A$ spectrum. Let for all sufficiently large $n$ the system of functions $f_n(x)$ be uniformly bounded and equicontinuous on this interval. Then the spectrum of $A$ is absolutely continuous on $[a,b\,]$, $\sigma'(x)\in C[a,b\,]$ and uniformly on this interval
$$
\lim\limits_{n\to\infty}\sigma_n(x)=\sigma(x)
$$
$$
\lim\limits_{n\to\infty} f_n(x)=\sigma'(x)
$$
\end{theorem}
-------------------------------------------------------------------------------

As an example consider here again Jacobi matrix associated with Hermite polynomials $\left(b_n=\sqrt{\frac{n+1}{2}}\,\right)$.
In this case for $|x|\le 2b_n$ one obtains
$$
f_n(x)=\frac{f_{J_n}(x)}{P^2_n(x)-\frac{\displaystyle b_{n-1}}{\displaystyle b_n}P_{n-1}(x)P_{n+1}(x)}=
\frac{b_n^2\,f_{J_n}(x)}{\sum\limits_{k=0}^n(b_k^2-b_{k-1}^2)P_k^2(x)}=
\frac{\sqrt{2(n+1)-x^2}}{\pi\sum\limits_{k=0}^n P_k^2(x)}=
\frac{h_n(x)}{g_n(x)}\,,
$$
where we denote for brevity the numerator and the denominator as $h_n(x)$ and $g_n(x)$. By simple transformations one can present the difference $f_n(x)-f_n(x')$ in the form
$$
f_n(x)-f_n(x')=f_n(x')\cdot\frac{g_n(x')-g_n(x)}{g_n(x)}+\frac{h_n(x)-h_n(x')}{g_n(x)}
$$
One has next
$$
\frac{h_n(x)-h_n(x')}{g_n(x)}=\frac{(x'-x)(x'+x)}{g_n(x)\left(\sqrt{2(n+1)-x^2}+\sqrt{2(n+1)-(x')^2}\right)}
$$
This expression tends to zero uniformly in $n$ as $x'\to x$.

For the second term one obtains
$$
f_n(x')\cdot\frac{g_n(x')-g_n(x)}{g_n(x)}=f_n(x')\cdot\frac{(x'-x)g_n'(\xi)}{g_n(x)}\,,
$$
where $\xi$ lies between $x$ and $x'$. From Из~(\ref{derivatives}) it follows that $|g_n'(x)|\le 4|x|\,g_n(x)$ for any $n$ and $x$. Hence using uniform boundedness of $f_n(x)$, monotonicity of $g_n(x)$ and letting $x'<x$, one obtains
$$
\left|f_n(x')\cdot\frac{g_n(x')-g_n(x)}{g_n(x)}\right|\le 4|\xi|\,|x'-x|\,f_n(x')\,\frac{g_n(\xi)}{g_n(x)}\le C|x'-x|
$$
Thus this expression also tends to zero uniformly as $|x'-x|\to 0$. Therefore the sequence $f_n(x)$ is equicontinuous on any interval of a real axis and the Theorem~(\ref{absolute_continuity_through_uniform_continuity}) gives
$$
\lim\limits_{n\to\infty} f_n(x)=\sigma'(x)
$$
or
$$
\lim\limits_{n\to\infty}\frac{\sqrt{2n}}{\pi\sum\limits_{k=0}^n P_k^2(x)}=\frac{1}{\sqrt{\pi}}\,e^{-x^2}\,,
$$
uniformly on any interval of a real axis. We can also write down this relation in terms of Hermite polynomials $H_n(x)$
$$
\lim\limits_{n\to\infty}\frac{\sqrt{2n}}{\pi\sum\limits_{k=0}^n\frac{\textstyle 1} {\textstyle k!\,2^k}\,H_k^2(x)}=\frac{1}{\sqrt{\pi}}\,e^{-x^2}\,,
$$
It is interesting to note that the functions
$$
y_k(x)=\frac{1}{\sqrt{\sqrt{\pi}\,k!\,2^k}}\,H_k(x)\,e^{-x^2/2}
$$
are normalized eigenfunctions of differential operator of harmonic oscillator $-\frac{d^2}{dx^2}+x^2$. This operator acting in
$L_2(\mathbb R)$ has a discrete spectrum with eigenvalues $\lambda_k=2k+1$. Hence the last limit equality yields
$$
\sum\limits_{k=0}^n y_k^2(x)\sim\frac{1}{\pi}\,\sqrt{\lambda_n}
$$
as $n\to\infty$.

Note that in the proof of equicontinuity of the functions $f_n(x)$ in this example we essentially use its boundedness.

Using the same reasoning one can obtain the following result:

\begin{theorem}
\label{condition_for_uniform_continuity}
Denote
$$
g_n(x)=P^2_n(x)-\frac{\displaystyle b_{n-1}}{\displaystyle b_n}P_{n-1}(x)P_{n+1}(x)
$$
Assume for all sufficiently large $n$ uniformly on the interval $[a,b\,]$ the estimate
$$
\frac{c_1}{b_n}<g_n(x)<\frac{c_2}{b_n}\,,\quad(0<c_1<c_2)\,,
$$
holds and
$$
|g_n'(x)|\le C\,g_n(x)\,,
$$
where the constant $C$ may depend on $x$ only. Then the sequence of functions $f_n(x)$ is uniformly bounded and equicontinuous on $[a,b\,]$.
\end{theorem}

-------------------------------------------------------------------------------

For the sequences $a_n$ and $b_n$ satisfying the conditions of the Theorem~(\ref{absolute_continuity_A_trough_a_n_b_n}), the first estimate of the Theorem~(\ref{condition_for_uniform_continuity}) is fulfilled. Make a hypothesis that the second estimate for the derivative $g_n'(x)$ is also fulfilled. By virtue of the Theorems~(\ref{absolute_continuity_and_convergence}) and ~(\ref{absolute_continuity_A_trough_a_n_b_n}) it would mean that $\sigma'(x)\in C(\mathbb R)$ and the sequence $f_n(x)$ converges to $\sigma'(x)$. The last fact would give the possibility to use the sequence $f_n(x)$ for a calculation of the spectral weight of the operator $A$.

\section{Conclusion.}

Note that the method of analysis of the absolutely continuous spectrum developed here can be applied not to the finite-difference operators only but to any self-adjoint operators in separable Hilbert space provided that one can find a convenient approximative sequence of operators $A_n$ with absolutely continuous spectrum strongly converging to the operator $A$ on a dense set.

Let us finish with the question which has remained open in this work: is the inverse statement to the Theorem~(\ref{absolute_continuity_and_convergence}) true? That is: do the uniform boundedness and equicontinuity of the functions $f_n(x)$ follow from the fact that the spectrum of $A$ is absolutely continuous on $[a,b\,]$ and $\sigma'(x)\in C[a,b\,]$?

\section{Acknowledgment.}

I would like to thank Grzegorz Swiderski for his very important remarks concerning the Theorems~(\ref{Hellinger}) and~(\ref{criteri_discretnosti}) and for the references on Chihara works and the works concerning Turan determinant.

\end{document}